%% file: locally-determined-colimits.tex
\documentclass[conference]{IEEEtran}
  \usepackage{fancyhdr}
  \usepackage{okgeneralmath}
  \usepackage{oksyntax}
  \usepackage{mathpartir}
  \usepackage{okcategory}
  \usepackage{oksemantics}
  \usepackage{oktheorem}
  \usepackage{float}
  \floatstyle{ruled}
  \usepackage{mathtools}

  \input{defs}


  \title{\vspace{-1cm}An absolute characterisation of\\locally determined $\omega$-colimits}

  \author{\IEEEauthorblockN{Ohad Kammar\\\texttt{<ohad.kammar@cl.cam.ac.uk>}}
    \IEEEauthorblockA{Programming, Logic, and Semantics
      Group\\University of Cambridge Computer Laboratory}}
\begin{document}
  \maketitle
  \begin{abstract}
    Characterising colimiting $\omega$-cocones of projection pairs in
    terms of least upper bounds of their embeddings and projections is
    important to the solution of recursive domain equations. We
    present a universal characterisation of this local property as
    $\omega$-cocontinuity of locally continuous functors. We present a
    straightforward proof using the enriched Yoneda embedding. The
    proof can be generalised to Cattani and Fiore's notion of locality
    for adjoint pairs.
  \end{abstract}

  \section*{Domains XI workshop contributed talk}
  In the category theoretic solution of recursive domain
  equations~\cite{smyth-plotkin:the-category-theoretic-solution-of-recursive-domain-equations},
  several technical results hinge upon the fact that the universality
  of $\omega$-cocones of projection pairs can be characterised
  \emph{locally} in terms of least upper bounds (lubs) of their
  embeddings and projections. To fix terminology and notation,
  consider an $O$-category $K$. Let $\PR K$ be the $O$-category
  consisting of \emph{projection pairs} $f \of A \to B$ given by $f =
  \pair{f^L \of A \to B}{f^R \of B \to A}$ where ${f^R \compose f^L =
    \id[A]}$ and ${f^L\compose f^R \leq \id[B]}$.
  \begin{definition*}[{\cite[Definition~8]{smyth-plotkin:the-category-theoretic-solution-of-recursive-domain-equations}}]
    We say that a cocone $\pair Cc$ for an $\omega$-chain of
    projection pairs is \emph{locally determined} when $\lub_{n \in
      \naturals}c_n^L\compose c_n^R = \id[C]$.

    When all colimiting $\omega$-cocones of projection pairs are
    locally determined, we say that the $O$-category \emph{has
      locally determined $\omega$-colimits of projection pairs}.
  \end{definition*}
  For example, the category $\wCPO$ of (not necessarily pointed)
  $\omega$-cpos and continuous functions has locally determined
  $\omega$-colimits.

  The importance of these cocones lies in the fact that every locally
  determined cocone is colimiting. As any locally continuous functor
  $F \of K \to L$ gives a continuous functor $\PR F \of \PR K \to \PR
  L$, given by $\PR F f \definedby \pair{Ff^L}{Ff^R}$, and locally
  determined $\omega$-cocones are preserved by these functors. Our
  contribution is to show the converse:
  \begin{theorem*}
    An $\omega$-colimiting cocone of projection pairs is locally
    determined if and only if it is preserved by every locally
    continuous functor.
  \end{theorem*}

  Let $\presheaf K$ be the $O$-category of $O$-presheaves, namely
  locally continuous functors and natural transformations from
  $\opposite K$ to $\wCPO$. Let $\yoneda \of K
  \to \presheaf K$ be the enriched Yoneda embedding $\yoneda x
  \definedby \wCPO(\placeholder, x)$. Then, following from general
  principles~\cite[Section~2.4]{kelly:basic-concepts-of-enriched-category-theory},
  $\yoneda$ is locally continuous and fully faithful.

  As is well-known, lubs and colimits in $O$-functor categories are
  given pointwise. The same argument shows that $\omega$-colimits of
  projection pairs are also given componentwise in $O$-functor
  categories. Therefore:
  \begin{proposition*}
    If $K$, $L$ are $O$-categories and $L$ has locally determined
    $\omega$-colimits of projection pairs, then so does the
    $O$-functor category $L^K$. In particular, every $O$-presheaf
    category $\presheaf K$ has locally determined $\omega$-colimits.
  \end{proposition*}

  We complete the proof of our theorem. Let $\pair Cc$ be any
  colimiting cocone that is preserved (in particular) by the locally
  continuous Yoneda embedding. As $\presheaf K$ has locally determined
  $\omega$-colimits:
  \[
  \yoneda\parent{\lub_{n}c_n^L\compose c_n^R} = \lub_n\yoneda (c_n^L) \compose \yoneda (c_n^R) = \yoneda (\id)
  \]
  By the faithfulness of the Yoneda embedding we deduce that $\pair
  Cc$ is locally determined.

  \begin{corollary*}
    An $O$-category has locally determined $\omega$-colimits of
    projection pairs if and only if every locally continuous functor
    from it yields an $\omega$-cocontinuous functor on projection
    pairs.
  \end{corollary*}

  Much of the theory of recursive domain equations generalises to
  \emph{adjoint pairs} $\pair{f^L}{f^R}$ where $f^L\compose f^R \leq
  \id$ and $\id \leq f^R \compose f^L$.  Cattani et
  al.~\cite{cattani-fiore-winskel:a-theory-of-recursive-domains-with-applications-to-concurrency,cattani-fiore:the-bicategory-theoretic-solution-of-recursive-domain-equations}
  generalised locally determined cocones as follows:
  \begin{definition*}[{cf.~\cite[Theorem~1.5]{cattani-fiore:the-bicategory-theoretic-solution-of-recursive-domain-equations}}]
    We say that a cocone $\pair Cc$ for an $\omega$-chain $\Delta$ of
    adjoint pairs is \emph{locally determined} when $\lub_{n \in
      \naturals}c_n^L\compose c_n^R = \id[C]$ and, for all $n \in
    \naturals$:
    \[
    \lub_{m \geq n} \Delta_{n \leq m}^R\compose \Delta_{m
      \geq n}^L = c_n^R\compose c_n^L
    \]
    When all colimiting $\omega$-cocones of adjoint pairs are
    locally determined, we say that the $O$-category \emph{has locally
      determined $\omega$-colimits of adjoint pairs}.
  \end{definition*}
  As $\wCPO$ has locally determined $\omega$-colimits of adjoint
  pairs, almost identical proofs show the following:
  \begin{theorem*}
    An $\omega$-colimiting cocone of adjoint pairs is locally
    determined if and only if it is preserved by every locally
    continuous functor.
  \end{theorem*}
  \begin{corollary*}
    An $O$-category has locally determined $\omega$-colimits of
    adjoint pairs if and only if every locally continuous functor from
    it yields an $\omega$-cocontinuous functor on adjoint pairs.
  \end{corollary*}

  \bibliographystyle{alpha}

\end{document}

%% file: defs.tex
\newcommand{\PR}[1]{{#1}_{\mathrm{PR}}}
\newcommand\lub{\bigvee}
\newcommand\presheaf[1]{\widehat{#1}}

%% file: locally-determined-colimits.bbl
\begin{thebibliography}{CFW98}

\bibitem[CF07]{cattani-fiore:the-bicategory-theoretic-solution-of-recursive-domain-equations}
Gian~Luca Cattani and Marcelo~P. Fiore.
\newblock The bicategory-theoretic solution of recursive domain equations.
\newblock {\em Electronic Notes in Theoretical Computer Science}, 172(0):203 --
  222, 2007.
\newblock Computation, Meaning, and Logic: Articles dedicated to Gordon
  Plotkin.

\bibitem[CFW98]{cattani-fiore-winskel:a-theory-of-recursive-domains-with-applications-to-concurrency}
G.L. Cattani, M.~Fiore, and G.~Winskel.
\newblock A theory of recursive domains with applications to concurrency.
\newblock In {\em Logic in Computer Science, 1998. Proceedings. Thirteenth
  Annual IEEE Symposium on}, pages 214--225, Jun 1998.

\bibitem[Kel82]{kelly:basic-concepts-of-enriched-category-theory}
Gregory~M. Kelly.
\newblock {\em Basic concepts of enriched category theory}.
\newblock Theory and Applications of Categories, 1982.
\newblock Reprinted in 2005.

\bibitem[SP82]{smyth-plotkin:the-category-theoretic-solution-of-recursive-domain-equations}
M.~Smyth and G.~Plotkin.
\newblock The category-theoretic solution of recursive domain equations.
\newblock {\em SIAM Journal on Computing}, 11(4):761--783, 1982.

\end{thebibliography}
